\documentclass{amsart}

\RequirePackage{hyperref}

\usepackage{amsmath,amsfonts}
\usepackage{graphicx}
\newtheorem{lem}{Lemma}
\newtheorem{conj}[lem]{Conjecture}
\newtheorem{prop}[lem]{Proposition}

\addtocounter{MaxMatrixCols}1
\begin{document}

%\begin{frontmatter}

%%  "Title of the Paper"
\title{Enumerative Geometry of the Mirror Quintic}%\protect\thanksref{T1}}
%\thankstext{T1}{???}

\author{Sheldon Katz}
\address{Department of Mathematics, University of Illinois, Urbana IL 61801}
\email{katzs@illinois.edu}
\author{David R. Morrison}
\address{Departments of Mathematics and Physics, University of California, Santa Barbara CA 93106}
\email{drm@math.ucsb.edu}

\begin{abstract}
We evaluate the enumerative invariants of low degree on the mirror
quintic threefold.
\end{abstract}

\maketitle

%\begin{keyword}[class=AMS]
%\kwd[Primary ]{}
%\kwd{}
%\kwd[; secondary ]{}
%\end{keyword}

%%  Upper case for every keyword
%\begin{keyword}
%\kwd{}
%\kwd{}
%\end{keyword}

%\tableofcontents

%\end{frontmatter}

%%  The body

Mirror symmetry burst upon the mathematical scene with the famous
computation made by Candelas, de la Ossa, Green, and Parkes 
\cite{cdgp} which in modern language proposed to count the number
of rational curves of fixed degree on the quintic threefold
using a technique from physics.  In current language
the ``counts'' are evaluations of Gromov--Witten invariants
or Gopakumar--Vafa invariants, and the technique for counting
these invariants has been expanded and extended in numerous ways.
From the physics point of view, the computation was made on
a closely related algebraic variety -- the mirror quintic.

In a recent physics paper \cite{HJKOV},
a new technique was proposed for explicitly evaluating the Gromov--Witten or
Gopakumar--Vafa invariants of the mirror quintic itself, not just the
quintic. Such explicit evaluations seem rather daunting, since
the answer will be a function of 101 variables.  One aspect of
\cite{HJKOV} is to use two variables only and arrive at a more
reasonable count.

Upon request of the authors of \cite{HJKOV}, the present authors worked out the
enumerative geometry of the mirror quintic.  We have made a conjecture
about the Mori cone which comes with a plausibility argument rather
than a proof.  However, indepedent of the truth of that conjecture,
we are able to explicitly evaluate the Gromov--Witten or
Gopakumar--Vafa invariants of low degree for the mirror quintic,
involving all 101 variables.

It gives us great pleasure to dedicate this paper to our mentor and friend Herb
Clemens.
Herb's  work on rational curves on Calabi--Yau threefolds
\cite{clemens1983double,clemensinf,MR934266} gave
an inspiration and a foundation to much of our own work, including this paper.
 
\section{The mirror quintic, its curves, and its divisors}\label{sec:curvesdivisors}
The standard description of the mirror quintic is as a quotient
of a quintic hypersurface in $\mathbb{P}^4$.  Let $t_1$, \dots, $t_5$
be homogeneous coordinates.  The standard equation is then
\begin{equation} t_1^5+t_2^5+t_3^5+t_4^5+t_5^5=5\psi t_1t_2t_3t_4t_5.\end{equation}
We are taking the quotient of $\mathbb{P}^4$
by the coordinate-wise action of 
\begin{equation} (\mathbb{Z}_5)^4 = \{ (\zeta_1,\zeta_2,\zeta_3,\zeta_4,\zeta_5) \in
(\mathbb{Z}_5)^5 \ | \ \prod_{j=1}^5 \zeta_j = 1\}.\end{equation}
The diagonal subgroup  acts trivially
on $\mathbb{P}^4$, so the effective
action is by $(\mathbb{Z}_5)^3$.

We will employ an alternate description (due to Batyrev \cite{Batyrev1993}) of
the mirror quintic as a singular complete intersection
in $\mathbb{P}^5$.  For this purpose we begin by describing
a basis for 
monomials invariant under the finite group action.  Our choice of
notation may appear awkward at first, but its utility will become
apparent later on.

Let
\begin{equation}
\begin{aligned}
u &= t_1t_2t_3t_4t_5 & v &= t_1^5 \\
w &= t_2^5 & x &= t_3^5\\
y &= t_4^5 & z &= t_5^5.
\end{aligned}
\end{equation}
Then the equations defining the mirror quintic can be written
\begin{align} \label{eq:first}
v + w + x + y + z &= 5\psi u \\ 
vwxyz &= u^5, \label{eq:second}
\end{align}
which describes the mirror quintic
 as a $(5,1)$ complete intersection in $\mathbb{P}^5$.  We denote this singular model of the mirror quintic
by $Y_\psi$, and refer to the coordinates appearing in \eqref{eq:first}
and \eqref{eq:second} as {\em Batyrev coordinates}.
In fact, this expresses the mirror quintic as a hypersurface (described
by \eqref{eq:first}) inside a singular toric variety (described by
\eqref{eq:second}).

There are singularities along lines defined by the simultaneous vanishing
of $u$ and two out of $\{v, w, x, y, z\}$ (as well as eq.~\eqref{eq:first}); 
there are ten such lines.
The transverse singularity at a general point of such a line is $A_4$.

There are more complicated singularities at the points defined by
the simultaneous vanishing of $u$ and three out of 
$\{v, w, x, y, z\}$ (and again imposing eq.~\eqref{eq:first}); 
there are ten such points.  The singularity at each such point
takes the local form $\mathbb{C}^3/(\mathbb{Z}_5)^2$, with non-isolated
singularities 
propagating along the coordinate axes in $\mathbb{C}^3$
having transverse singularity of the form $\mathbb{C}^2/\mathbb{Z}_5 = A_4$.

We also observe the feature of five $\mathbb{P}^2$'s contained within
the mirror quintic:  each is given by the simultaneous vanishing of $u$
and one of the coordinates $\{v, w, x, y, z\}$, as well as
eq.~\eqref{eq:first}.

\begin{figure}
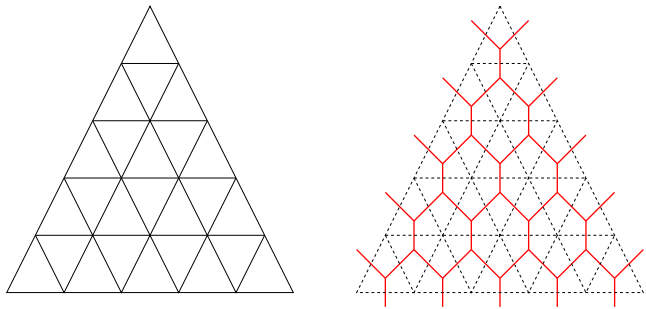

\begin{center}
\includegraphics[scale=0.3]{5trianglea.mps}
\qquad
\includegraphics[scale=0.3]{5triangle3b.mps}
\end{center}
\caption{A triangulation specifying a crepant resolution 
of $\mathbb{C}^3/(\mathbb{Z}_5)^2$, and
the corresponding dual graph.}\label{fig:triangulation}
\end{figure}

A concrete procedure for blowing up the singularities was given
in \cite[Appendix B]{mirrorguide}, but we will use a different
birational model described in \cite[Section
4]{geomaspects}.\footnote{Note that our use of the term ``the mirror
quintic'' is not quite correct -- we must specify the birational model,
and different models will have different K\"ahler cones.  In this
paper, we use the model corresponding to the resolution of singularities
specified in \cite[Section 4]{geomaspects}.}  The
singularity $\mathbb{C}^3/(\mathbb{Z}_5)^2$ can be described
torically, and we display the toric data for the resolution we
use in the left half of Figure~\ref{fig:triangulation}.  We will denote
the resulting Calabi-Yau threefold by $X_\psi$.  The fully-resolved 
mirror quintic family is parametrized by $z=\psi^{-5}$, and is smooth for 
$z\in\mathbb{P}^1-\{0,1,\infty\}$. The large radius point is given by $z=0$, 
while $z=1$ is the conifold point and  $z=\infty$ is the orbifold point.  There is a 
birational contraction map $\rho:X_\psi\to Y_\psi$.

The toric data is ``dual'' to a more geometric description involving
curves and surfaces, and we display the corresponding dual graph
in the right half of Figure~\ref{fig:triangulation}.  There are six
complete divisors within the resolution, represented by hexagons
in the dual graph.  There are also four incomplete divisors along
each edge which represent components of the
resolutions of the one-dimensional singular
loci, and three incomplete divisors at the vertices of the figure
which represent the $\mathbb{P}^2$'s identified above.

There are also $30$ compact curves within the resolution, represented
by segments in the dual graph which have both ends meeting a divisor.
We can also see portions of other curves which we will describe later.

\begin{figure}
\begin{center}
\includegraphics[scale=0.5]{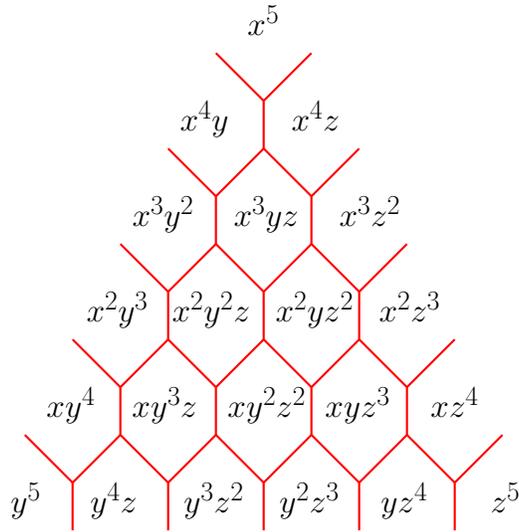}
\end{center}
\caption{Labeling the divisors with monomials.}
\label{fig:triangulation-labeled}
\end{figure}

In order to keep track of all of the divisors and curves,  we
introduce the following notation.  Each point with stabilizer
$\mathbb{C}^5/(\mathbb{Z}_5)^2$ is associated with three of the
variables $\{v, w, x, y, z\}$ and we will use quintic monomials
in those three variables
to label the corresponding
divisors.
This is illustrated in 
Figure~\ref{fig:triangulation-labeled} for the variables $x$, $y$, $z$.
For a given monomial $\mathbf{m}$, we call the corresponding divisor
$D_{\mathbf{m}}$.

Note that every monomial which involves at most three of 
%the variables
$\{v,w,x,y,z\}$ corresponds to a divisor, and monomials involving fewer
than three of 
%the 
those
variables will show up in more than one of the 
toric diagrams.  There are 105 divisors of this form.

Each compact curve in our diagram is the intersection of precisely
two compact divisors $D_{\mathbf{m}}$ and $D_{\mathbf{n}}$ and
we label the curve by\footnote{Our convention is to only use the
notation $\gamma_{\mathbf{m},\mathbf{n}}$ in case the two divisors
$D_{\mathbf{m}}$ and $D_{\mathbf{n}}$ are both compact in the inverse
image of some $\mathbb{C}^5/(\mathbb{Z}_5)^2$ point.} $\gamma_{\mathbf{m},\mathbf{n}}$.  In addition
for each variable $\mathbf{s}$ from our set of variables 
$\{v, w, x, y, z\}$, we let $\ell_{\mathbf{s}}$ denote a line in
$D_{\mathbf{s}^5} \cong \mathbb{P}^2$.  (Note that $\ell_{\mathbf{s}}$ is also represented by
the curve $D_{\mathbf{s}^5} \cap D_{\mathbf{s}^4\mathbf{t}}$ 
for any variable $\mathbf{t}$ distinct from $\mathbf{s}$.)
Finally, given two variables $\mathbf{s}$ and $\mathbf{t}$ from
the set $\{v, w, x, y, z\}$, we let $\sigma_{\mathbf{s},\mathbf{t}}$
be the intersection of
$D_{\mathbf{s}^3\mathbf{t}^2}$ and
$D_{\mathbf{s}^2\mathbf{t}^3}$.
We will conjecture below that
the 315 curve classes of the form $[\gamma_{\mathbf{m},\mathbf{n}}]$,
$[\ell_{\mathbf{s}}]$, and $[\sigma_{\mathbf{s},\mathbf{t}}]$
 are precisely the generators of all
extremal rays and collectively generate
the Mori cone of the mirror quintic.

There are a few more curves which should be discussed.  For each pair of
variables such as $\{y,z\}$, resolving the corresponding one-dimensional
 singular locus
produces four divisors which in the example are
$D_{y^4z}$, $D_{y^3z^2}$, $D_{y^2z^3}$, and $D_{yz^4}$.  Each divisor
appears over three of the $\mathbb{C}^3/(\mathbb{Z}_5)^2$ points
(corresponding to the triples of variables which contain 
$\{y,z\}$) and they are labeled the same way in each instance.
These divisors are all ruled (by the exceptional curves resolving
the non-isolated singularities), and we label the fiber $\varphi$ of the ruling
with the same monomial as that of the divisor (e.g., the ruling on $D_{y^4z}$
is
$\varphi_{y^4z}$). 
There are also  curves 
%$\sigma$
of intersection of adjacent pairs
of such divisors, labeled by a pair of monomials, e.g., the intersection
of $D_{y^4z}$ and $D_{y^3z^2}$.
The only extremal curves among those are the lines 
$\ell_{\mathbf{s}}$, $\ell_{\mathbf{t}}$ and the ``middle'' sections
$\sigma_{\mathbf{s},\mathbf{t}}$.

\begin{figure}
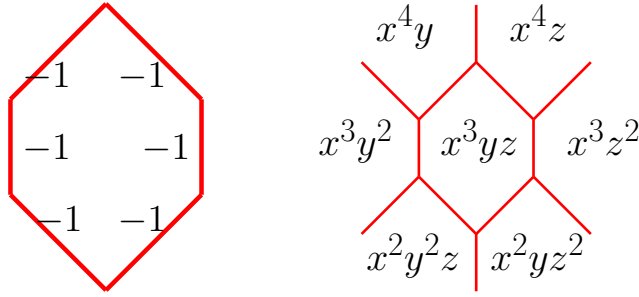

\begin{center}
\includegraphics[scale=1.0]{hex.mps}
\qquad
\qquad
\includegraphics[scale=0.6]{hex-lbl.mps}
\end{center}
\caption{A general $dP_3$, and a particular $dP_3$.}\label{fig:one-divisor}
\end{figure}

We determine the relations among the various curve classes by 
analyzing the structure of the divisors.  Each compact divisor
of the form $D_{\mathbf{m}}$ where the monomial $\mathbf{m}$ involves
three variables (which is represented by a hexagon in the dual graph) is
isomorphic to $dP_3$, i.e., the blowup of $\mathbb{P}^2$ in three
non-collinear points.  Each of the six toric curves on the divisor
has self-intersection $-1$, as illustrated in the left half of
Figure~\ref{fig:one-divisor}.  In the right half of that same
figure, we have chosen a particular example so as to provide notation
for the relations on curve classes that we are about to describe.

\begin{figure}
\begin{center}
\includegraphics[scale=0.5]{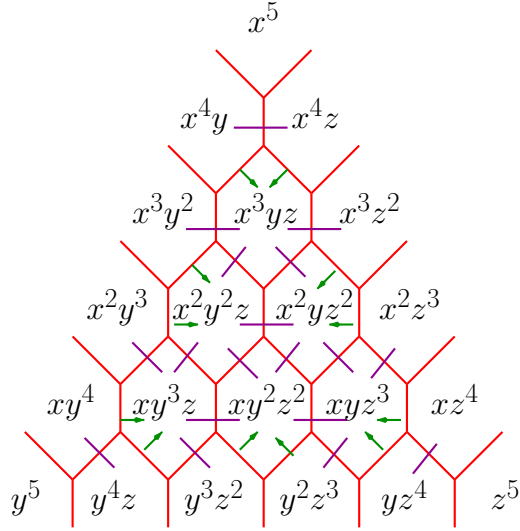}
\end{center}
\caption{Basis for compact curves.}
\label{fig:basis}
\end{figure}

A $dP_3$ has three rulings, each having two singular fibers.  Each singular
fiber in one of these rulings is the union of two adjacent $-1$ curves.
For example, the ruling whose general fibers run from southwest to 
northeast in the right side of Figure~\ref{fig:one-divisor} has
special fibers 
$\gamma_{x^4y,x^3yz} \cup \gamma_{x^3y^2,x^3yz}$ and
$\gamma_{x^3yz,x^3z^2} \cup \gamma_{x^3yz,x^2yz^2}$.
This leads to a relation among the curve classes which we denote by
$R^{xy,z}_{x^3yz}$, namely
\begin{equation}
R^{xy,z}_{x^3yz} = [\gamma_{x^4y,x^3yz}] + [\gamma_{x^3y^2,x^3yz}]
  - [\gamma_{x^3yz,x^3z^2}] - [\gamma_{x^3yz,x^2yz^2}] .
\end{equation}
The subscript on $R$ denotes the divisor which generated the relation.
The superscript identifies the ruling in the following way:  the parameter
space for the ruling proceeds from the $xy$ edge of the toric diagram to
the $z$ vertex of the toric diagram.  In this way, we can list all such
relations in the $xy,z$ direction, as follows:
\begin{equation} \label{eq:relgamma}
\begin{aligned}
R^{xy,z}_{x^3yz} &= [\gamma_{x^4y,x^3yz}] + [\gamma_{x^3y^2,x^3yz}]
  - [\gamma_{x^3yz,x^3z^2}] - [\gamma_{x^3yz,x^2yz^2}] \\
R^{xy,z}_{x^2y^2z} &= [\gamma_{x^3y^2,x^2y^2z}] + [\gamma_{x^23y^3,x^2y^2z}]
  - [\gamma_{x^1y^2z,x^2yz^2}] - [\gamma_{x^2y^2z,xy^2z^2}] \\
R^{xy,z}_{xy^3z} &= [\gamma_{x^2y^3,xy^3z}] + [\gamma_{xy^4,xy^3z}]
  - [\gamma_{xy^3z,xy^2z^2}] - [\gamma_{xy^3z,y^3z^2}] \\
R^{xy,z}_{x^2yz^2} &= [\gamma_{x^3yz,x^2yz^2}] + [\gamma_{x^2y^2z,x^2yz^2}]
  - [\gamma_{x^2yz^2,x^2z^3}] - [\gamma_{x^2yz^2,xyz^3}] \\
R^{xy,z}_{xy^2z^2} &= [\gamma_{x^2y^2z,xy^2z^2}] + [\gamma_{xy^3z,xy^2z^2}]
  - [\gamma_{xy^2z^2,xyz^3}] - [\gamma_{xy^2z^2,y^2z^3}] \\
R^{xy,z}_{xyz^3} &= [\gamma_{x^2yz^2,xyz^3}] + [\gamma_{xy^2z^2,xyz^3}]
  - [\gamma_{xyz^3,xz^4}] - [\gamma_{xyz^3,yz^4}] \\
\end{aligned}
\end{equation}
Similarly, we can relate the 
fiber in the ruling on one of the ``edge'' divisors to a degenerate
fiber within the toric diagram.  The ones in the same $xy,z$ direction are:
\begin{equation} \label{eq:relphi}
\begin{aligned}
R^{xy,z}_{x^4y} &= 
  [\varphi_{x^4y}] - [\gamma_{x^4y,x^4z}] - [\gamma_{x^4y,x^3yz}] \\
R^{xy,z}_{x^3y^2} &= 
  [\varphi_{x^3y^2}] - [\gamma_{x^3y^2,x^3yz}] - [\gamma_{x^3y^2,x^2y^2z}] \\
R^{xy,z}_{x^2y^3} &= 
  [\varphi_{x^2y^3}] - [\gamma_{x^2y^3,x^2y^2z}] - [\gamma_{x^2y^3,xy^3z}] \\
R^{xy,z}_{xy^4} &= 
  [\varphi_{xy^4}] - [\gamma_{xy^4,xy^3z}] - [\gamma_{xy^4,y^4z}] \\
\end{aligned}
\end{equation}
Other relations of the same type as in eqs.~\eqref{eq:relgamma} and
\eqref{eq:relphi} can be generated by substituting any three variables
for $x$, $y$, $z$.  This gives $30$ versions of eq.~\eqref{eq:relgamma}.

The compact divisor $D_{x^3yz}$ also allows us to exhibit a relation among
relations, since the three specified relations among six $-1$ curves are not 
linearly independent.  The ``syzygy'' is easily seen to be
\begin{equation}
R^{xy,z}_{x^3yz} + R^{yz,x}_{x^3yz} + R^{zx,y}_{x^3yz} = 0.
\end{equation}
A similar syzygy exists for each compact divisor.
We thus conclude that of the 180 relations of type \eqref{eq:relgamma}
among the 300 $\gamma$ curves, only 120 are linearly independent.

We can implement these relations as illustrated in
Figure~\ref{fig:basis}.  In that figure, we have selected a subset of the 
compact curves (marked with a purple dash) which constitute $4$ curves
on each compact divisor.  In addition, for the unmarked compact curves,
we have used a green arrow to designate which divisor's relations
should be used to eliminate that curve from the spanning set.
Each compact divisor has two green arrows within it, indicating that two
of the relations belonging to that divisor are used in this elimination
process.  Those relations are linearly independent.
 
The divisors on the edges, such as $D_{x^4y}$ also give rise to
relations among the $\gamma$ curves, since three expressions for
$[\varphi_{x^4y}]$ are derived from relations of the type in
eq.~\eqref{eq:relphi}.  The ones for $\varphi_{x^4y}$ are
\begin{equation}
\begin{aligned}
R^{xy,v}_{x^4y} &= 
  [\varphi_{x^4y}] - [\gamma_{x^4y,x^4v}] - [\gamma_{x^4y,x^3yv}] \\
R^{xy,w}_{x^4y} &= 
  [\varphi_{x^4y}] - [\gamma_{x^4y,x^4w}] - [\gamma_{x^4y,x^3yw}] \\
R^{xy,z}_{x^4y} &= 
  [\varphi_{x^4y}] - [\gamma_{x^4y,x^4z}] - [\gamma_{x^4y,x^3yz}] \\
\end{aligned}
\end{equation}
and by subtraction we obtain two independent relations among the $\gamma$
curves.
Since this can be done for each of 4 divisors along each of 10 edges, there
are 80 relations in total of this type.  These relations are independent.

Thus, we find a total of 200 linear relations among the 300 compact 
$\gamma$ curves, leaving 100 independent classes.

\begin{figure}
\begin{center}
\includegraphics[scale=1.0]{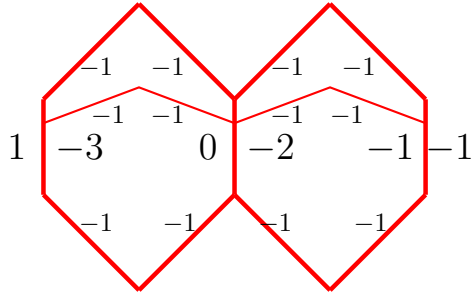}
\end{center}
\caption{Two adjacent components of the resolution of the $A_4$ locus.}\label{fig:two-divisors}
\end{figure}

We now consider the divisors on the edges of the toric diagram, which form the
resolution of the $A_4$ loci.  We take as an example the sequence of
divisors $D_{y^4z}$, $D_{y^3z^2}$, $D_{y^2z^3}$, $D_{yz^4}$.  The
first two divisors in any such sequence are illustrated in
Figure~\ref{fig:two-divisors}, in other words, the divisor on the
left is $D_{y^4z}$ and the divisor on the right is $D_{y^3z^2}$.  To
the left of $D_{y^4z}$ is $D_{y^5}$ and to the right of $D_{y^3z^2}$
is $D_{y^2z^3}$.  (There is also a third singular fiber in the ruling on each 
divisor; these are indicated with narrow lines.)

As we mentioned earlier, the divisors $D_{y^5}$ and $D_{y^4z}$ meet
on $\ell_y$ which has self-intersection $1$ on $D_{y^5}$.  By the
adjunction formula, the self-intersection of this curve on $D_{y^4z}$
(where it is a section of the ruling)
must be $-3$. It follows that $D_{y^4z}$ is a blowup of the Hirzebruch
surface $\mathbb{F}_3$, and since there are three singular fibers, the
disjoint section $D_{y^4z}\cap D_{y^3z^2}$ must have self-intersection $0$
on $D_{y^4z}$.  In fact, we can say more:  by the description as
a blowup of a Hirzebruch surface\footnote{If we blowup the Hirzebruch surface
$\mathbb{F}_n$ at $k$ points which are contained in a section $\overline{\sigma}$
of the $\mathbb{P}^1$-fibration (with fiber $\overline{\varphi}$) having
$(\overline{\sigma})^2_{\mathbb{F}_n}=n\ge0$, then the blowup contains the
proper transform $\varphi$ of $\overline{\varphi}$, the 
proper transform $\sigma_\infty$ of the section $\overline{\sigma_\infty}$
with $(\overline{\sigma_\infty})^2_{\mathbb{F}_n}=-n\le0$, the proper transform
$\sigma$ of $\overline{\sigma}$, and the exceptional divisors $e_1$, \dots, $e_k$.
The basic homology relation
$[\overline{\sigma}]=[\overline{\sigma_\infty}]+n[\overline{\varphi}]$ on $\mathbb{F}_n$
pulls back to a homology relation
$[\sigma]+[e_1]+\cdots+[e_k]=[\sigma_\infty]+n[\varphi]$
on the blowup.  Note that $[\varphi]-[e_j]$ is represented by an effective $(-1)$-curve
$e_j'$, and that $[e_j]+[e_j']$ is homologous to $[\varphi]$.}
\begin{equation}
 [D_{y^4z}\cap D_{y^3z^2}] = [\ell_y] + 3[\varphi_{y^4z}] 
- [\gamma_{vy^3z,y^4z}]
- [\gamma_{wy^3z,y^4z}]
- [\gamma_{xy^3z,y^4z}],
\end{equation}
which can also be written
\begin{equation}
 [D_{y^4z}\cap D_{y^3z^2}] = 
[\ell_y] 
+ [\gamma_{vy^4,y^4z}]
+ [\gamma_{wy^4,y^4z}]
+ [\gamma_{xy^4,y^4z}].
\label{eq:sectioninMori}
\end{equation}

Using adjunction again, the self-intersection of $D_{y^4z}\cap D_{y^3z^2}$
on $D_{y^3z^2}$ must be $-2$, which tells us that $D_{y^3z^2}$ is the
blowup of a Hirzebruch surface $\mathbb{F}_2$ in three points.  The disjoint
section $\sigma_{y,z}= D_{y^3z^2}\cap D_{y^2z^3}$ 
must thus have self-intersection
$-1$ on $D_{y^3z^2}$, and again we get 
a relation among curve classes:
\begin{equation}
 [\sigma_{y,z}] = [D_{y^4z}\cap D_{y^3z^2}] + 2[\varphi_{y^3z^2}] 
- [\gamma_{vy^2z^2,y^3z^2}]
- [\gamma_{wy^2z^2,y^3z^2}]
- [\gamma_{xy^2z^2,y^3z^2}],
\end{equation}
which can also be written
\begin{equation}
 [\sigma_{y,z}] =  [D_{y^4z}\cap D_{y^3z^2}]
+ [\gamma_{vy^3z,y^3z^2}]
+ [\gamma_{wy^3z,y^3z^2}]
- [\gamma_{xy^2z^2,y^3z^2}]
\end{equation}
(a little less naturally this time, since we must single out one of the
three singular fibers).
Combining these two gives
\begin{equation}
 [\sigma_{y,z}] =  
[\ell_y] 
+ [\gamma_{vy^4,y^4z}]
+ [\gamma_{wy^4,y^4z}]
+ [\gamma_{xy^4,y^4z}]
+ [\gamma_{vy^3z,y^3z^2}]
+ [\gamma_{wy^3z,y^3z^2}]
- [\gamma_{xy^2z^2,y^3z^2}].
\label{eq:sigma}
\end{equation}

%On the other hand, this same curve class $[\sigma_{y^4z,y^3z^2}]$
%can be expressed in terms of $[\ell_z]$ and the $\gamma$
%curves by
%starting at the ``$z$ end'' of the edge, producing the relation
%\begin{equation}
% [\sigma_{y^4z,y^3z^2}] =  
%[\ell_z] 
%+ [\gamma_{vz^4,yz^4}]
%+ [\gamma_{wz^4,yz^4}]
%+ [\gamma_{xz^4,yz^4}].
%+ [\gamma_{vyz^3,y^2z^3}]
%+ [\gamma_{wyz^3,y^2z^3}]
%- [\gamma_{xy^2z^2,y^2z^3}].
%\end{equation}
%Taking the difference of these two expresses the difference of 
%$[\ell_z]$ and $[\ell_y]$ in terms of the $\gamma$ curves:
%\begin{equation}
%\begin{aligned}
%{} [\ell_z]
%-[\ell_y] = 
%&- [\gamma_{vz^4,yz^4}]
%- [\gamma_{wz^4,yz^4}]
%- [\gamma_{xz^4,yz^4}].
%- [\gamma_{vyz^3,y^2z^3}]
%- [\gamma_{wyz^3,y^2z^3}]
%+ [\gamma_{xy^2z^2,y^2z^3}] \\
%&+ [\gamma_{vy^4,y^4z}]
%+ [\gamma_{wy^4,y^4z}]
%+ [\gamma_{xy^4,y^4z}].
%+ [\gamma_{vy^3z,y^3z^2}]
%+ [\gamma_{wy^3z,y^3z^2}]
%- [\gamma_{xy^2z^2,y^3z^2}].
%\end{aligned}
%\end{equation}

We thus get 20 additional relations among our 315 curves.  However, at
most
14 of the 20 relations can be linearly independent, which can be seen
from Batyrev's calculation of $b_4$ of the quintic mirror
\cite{Batyrev1993}.

The conclusion (from Batyrev's calculation) is that there are  precisely
101 independent curve classes generated by our explicit set of 315 curves.

\bigskip

We now consider relations among the divisors.  First, for the divisors
whose monomials are chosen from $\{x,y,z\}$, we can determine some
linear combinations which have intersection number $0$ with all 30
compact curves lying over the corresponding singular point.  There
are three linearly independent combinations, which we label as 
$D^{x,yz}$,
$D^{y,xz}$, and $D^{z,xy}$, and which are described by means of
coefficient arrays whose shape matches that of 
Figure~\ref{fig:triangulation-labeled},
as follows:
\[ 
D^{y,xz} = \begin{matrix}
&&&&& 0 &&&&& \\
&&&& 1 && 0 &&&& \\
&&& 2 && 1 && 0 &&& \\
&& 3 && 2 && 1 && 0 &&\\
& 4 && 3 && 2 && 1 && 0 &\\
5 && 4 && 3 && 2 && 1 && 0\\
\end{matrix}
, 
\]
%\qquad
\[
D^{z,xy} = \begin{matrix}
&&&&& 0 &&&&& \\
&&&& 0 && 1 &&&& \\
&&& 0 && 1 && 2 &&& \\
&& 0 && 1 && 2 && 3 &&\\
& 0 && 1 && 2 && 3 && 4 &\\
0 && 1 && 2 && 3 && 4 && 5\\
\end{matrix}
, 
\]
\[
D^{x,yz} = \begin{matrix}
&&&&& 5 &&&&& \\
&&&& 4 && 4 &&&& \\
&&& 3 && 3 && 3 &&& \\
&& 2 && 2 && 2 && 2 &&\\
& 1 && 1 && 1 && 1 && 1 &\\
0 && 0 && 0 && 0 && 0 && 0\\
\end{matrix}
.
\]
Note that $\ell_x$ meets $D_{x^5}$ with intersection number $-3$, 
and also meets $D_{x^4y}$ and $D_{x^4z}$, each with
intersection number $1$.  It follows that $D^{x,yz}$ meets $\ell_x$ with intersection number $-7$, meets $\ell_y$ once, and meets $\ell_z$ once.

This calculation suggests how to find a divisor which has intersection number $0$
with all curves $\gamma_{\mathbf{m_1},\mathbf{m_2}}$.  Let
\begin{equation}
D^{x} = \sum_{\substack{\text{monomials } \mathbf{m} \text{ with}\\
\text{at most three variables}}} 
\nu_x(\mathbf{m})\, D_{\mathbf{m}},
\end{equation}
where $\nu_x(\mathbf{m})$ denotes the exponent of the variable $x$ in the monomial $\mathbf{m}$.
Since the coefficients of $D^x$ in any triangle containing $x$ have the
same pattern as $D^{x,yz}$, the intersection with all of the compact
curves over such a point is $0$.  And since the coefficients of $D^x$
are all zero in any triangle not containing $x$, the intersection
number with each of the compact curves over such a point is also $0$.

Also, $D^x$ meets $\ell_x$ at the components $D_{x^5}$ and $D_{x^4*}$,
so the total intersection number is 
\[ D^x\cdot \ell_x =5\cdot(-3) + 4\cdot4\cdot1 = 1.\]
On the other hand, for any other variable such as $y$, the only component
of $D^x$ with a nonzero coefficient which meets $\ell_y$ is $D_{x^4y}$.
Thus, $D^x\cdot \ell_y=1$ as well (and similarly for $\ell_v$, $\ell_w$,
and $\ell_z$).

It follows that $D^x\cdot\gamma=0$, $D^x\cdot \ell=1$ for all $\gamma$
curves and lines $\ell$, respectively.

Now, the same thing happens for $D^v$, $D^w$, $D^y$, and $D^z$, so the
difference of any two is numerically trivial.  This produces four linear
relations among the 105 divisors $D_{\mathbf{m}}$, giving a space of
dimension 101.

Note that $D^x$ (and the others) can be identified with a hyperplane section
$H$ of the original singular model of the mirror quintic.

\bigskip

We would like to write down one more divisor, which has intersection number
$1$ with all $\gamma$ curves, and intersection number $0$ with all lines
$\ell$.  To this end, we consider the following coefficient array
\[ 
\begin{matrix}
&&&&& 16 &&&&& \\
&&&& 12 && 12 &&&& \\
&&& 10 && 9 && 10 &&& \\
&& 10 && 8 && 8 && 10 &&\\
& 12 && 9 && 8 && 9 && 12 &\\
16 && 12 && 10 && 10 && 12 && 16\\
\end{matrix}
\]
The corresponding divisor
 has intersection number $1$ with each compact $\gamma_{\mathbf{m}}$
where $\mathbf{m}$ involves the variables $\{x,y,z\}$.  
This is verified by checking total intersection number in various
sub-diamonds of the coefficient array, such as $16+9-12-12$ which
corresponds to intersections with $\gamma_{x^4y,x^4z}$.

To make a global
version of this, we define a function $a$ on monomials as follows.
The nonzero exponents in the monomial $\mathbf{m}$ determine a partition
of $5$ into at most three elements, and the value of the function
only depends on the partition, as follows:

\noindent
\begin{center}
\begin{tabular}{c|c|c|c|c|c}
Partition & $5$ & $4+1$ & $3+2$ & $3+1+1$ & $2+2+1$ \\ \hline
$a(\mathbf{m})$ & $16$ & $12$ & $10$ & $9$ & $8$
\end{tabular}
\end{center}

\noindent
Then the divisor 
\[ \mathbb{D} = 
\sum_{\substack{\text{monomials } \mathbf{m} \text{ with}\\
\text{at most three variables}}} 
a(\mathbf{m})\, D_{\mathbf{m}}
\] has the property that $\mathbb{D}\cdot \gamma=1$ for all $\gamma$.
Note that $\ell_x$ meets $D_{x^5}$ of coefficient $16$ with intersection
number $-3$, and meets 4 divisors of the form $D_{x^4*}$ each of coefficient
$12$ and with intersection number $1$.  Thus,
\[ \mathbb{D}\cdot \ell_x = 16\cdot (-3) + 4 \cdot 12 \cdot 1 = 0\]
and similarly for any of the $\mathbb{D}\cdot\ell$.

\bigskip
Having identified the curves $\ell_{\mathbf{t}}$ and $\sigma_{\mathbf{s},\mathbf{t}}$ for $\mathbf{s},\mathbf{t}\in\{w,v,x,y,z\}$, and the curves $\gamma_{\mathbf{m_1},\mathbf{m_2}}$ above,
we conjecture that these generate the Mori cone of numerical equivalence classes of effective curves on $X_\psi$.

\bigskip\noindent
\begin{conj}
The Mori cone $M$ of $X_\psi$ is generated by the classes of the curves $\ell_{\mathbf{t}}$, $\sigma_{\mathbf{s},\mathbf{t}}$, and  $\gamma_{\mathbf{m_1},\mathbf{m_2}}$.
\label{conj:main}
\end{conj}

\medskip
\emph{In the rest of this paper, we assume this conjecture whenever necessary.}

\bigskip
If $C$ is a curve on $X_\psi$, we define the \emph{degree} of $C$ to be the degree of its image $\overline{C}:=\rho(C)\subset Y_\psi\subset\mathbb{P}^5$. We will show presently in Lemma~\ref{lem:deg0} that if $C\subset X_\psi$ has degree zero, then its class $[C]$ is in the cone generated by the classes of the curves $\gamma_{\mathbf{m_1},\mathbf{m_2}}$.  Furthermore, in Section~\ref{sec:enum} we will show that if $C\subset X_\psi$ has degree 1, then its class $[C]$ is in the cone generated by the classes of the curves $\ell_{\mathbf{t}}$, $\sigma_{\mathbf{s},\mathbf{t}}$, and  $\gamma_{\mathbf{m_1},\mathbf{m_2}}$.  In other words, the conjecture is true for curves of degree at most 1.

\begin{lem}
The class of any curve in an exceptional divisor for the birational contraction $\rho$ is in the cone generated by the classes of the curves $\ell_{\mathbf{t}}$, $\sigma_{\mathbf{s},\mathbf{t}}$, and  $\gamma_{\mathbf{m_1},\mathbf{m_2}}$.
\label{lem:deg0}
\end{lem}

In particular, the conjecture is true for curves of degree~0.

\bigskip\noindent
\emph{Proof.} We give separate arguments for the exceptional divisors which contract to a point and to a curve.  

The exceptional divisors contracting to a point are all $dP_3$s.  It is well-known that the Mori cone of $dP_3$ is generated by the classes of its six exceptional curves of the first kind.  But these exceptional curves are precisely the $\gamma$ curves contained in that $dP_3$.

The exceptional divisors contracting to a curve are the components of one of the $A_4$ resolutions. By symmetry we need only analyze the curves in $D_{y^4z}$ and $D_{y^3z^2}$.

Now $D_{y^4z}$ is identified with the blow up of the Hirzebruch surface
$\mathbb{F}_3$ at three general points, the exceptional curves being
identified with the curves $\gamma_{y^4z,y^3z\mathbf{s}}$ for
$\mathbf{s}\in\{v,w,x\}$.  Using the blowup description, the Mori cone
is generated by the $-3$ section, the proper transform of any
$+3$ section containing the points being blown up, the proper
transform of the fibers, and the exceptional curves.  We need only
check each curve in turn.  The $-3$ section is identified with the
curve $\ell_y$.  We have already noted that the exceptional curves are
all $\gamma$ curves.  The proper transform of the fiber which meets
$\gamma_{y^4z,y^3z\mathbf{s}}$ is just $\gamma_{y^4z,y^4\mathbf{s}}$,
also a $\gamma$ curve.  The proper transform of a particular $+3$ section is
identified with $D_{y^4z}\cap D_{y^3z^2}$, which is in the cone
generated by $\ell_y$ and $\gamma$ curves by (\ref{eq:sectioninMori}).

Next, we observe that $D_{y^3z^2}$ is similarly identified with the blow up of the Hirzebruch surface
$\mathbb{F}_2$ at three points on a $+2$ section, the exceptional curves being
identified with the curves $\gamma_{y^3z^2,y^2z^2\mathbf{s}}$ for
$\mathbf{s}\in\{v,w,x\}$.  Using the blowup description, the Mori cone
is generated by the $-2$ section, the proper transform of the unique
$+2$ section containing the points being blown up, the proper
transform of the fibers, and the exceptional curves.  We need only
check each curve in turn.  The $-2$ section is identified with the
curve $D_{y^4z}\cap D_{y^3z^2}$ already considered above.  
We have already noted that the exceptional curves are
all $\gamma$ curves.  The proper transform of the fiber which meets
$\gamma_{y^4z,y^3z\mathbf{s}}$ is just $\gamma_{y^4z,y^4\mathbf{s}}$,
also a $\gamma$ curve.  The proper transform of the $+2$ section is
identified with $\sigma_{y,z}$.

\bigskip
We can give a plausibility argument in favor of our conjecture for arbitrary degree as follows.
%I would like to argue in favor of the Mori cone statement as follows.
Let $C$ be any irreducible curve on the mirror quintic, and let 
$\overline{C}$ be its image on the singular model as above.  If  $\overline{C}$
is contained in the singular locus, then $C$ must lie in one or more
of the exceptional divisors of the blowup.  This case is handled by Lemma~\ref{lem:deg0}. If $\overline{C}$ is not contained in
the singular locus then it meets the singular locus in finitely many points.
We want to claim -- and this is the gap in the argument -- that $\overline{C}$
can be deformed away from the singular locus.  The family $C_t$, when
lifted to the smooth model, will have a limit $C_0$ consisting of $C$
together with some curves which are contained in exceptional divisors.
As mentioned earlier, those latter curves are in the cone generated
by known curves.

\bigskip
We next let $\mathbf{H}=\mathbb{Z}^2=\mathbb{Z}\cdot \ell\oplus\mathbb{Z}\cdot\gamma$ and project curve classes to $H$ via
\begin{equation}\label{eq:proj2}
\pi:M\to \mathbf{H},\qquad \pi(C)=\left(C\cdot D^x\right)\ell+\left(C\cdot\mathbb{D}\right)\gamma.
\end{equation}
The notation for the basis for $\mathbf{H}$ has been chosen so that $\pi(\ell_{\mathbf{t}})=\ell$ and $\pi(\gamma_{\mathbf{m_1},\mathbf{m_2}})=\gamma$.  We also observe that if $\pi(\beta)=m\ell+n\gamma$, then $\rho_*(\beta)$ is the class of a 
degree $m$ curve in $Y_\psi\subset\mathbb{P}^5$.

\smallskip
We compute and record the following curve classes, which are computed from (\ref{eq:proj2}) and the definitions of
$D^x$ and $\mathbb{D}$.

\begin{eqnarray}\label{eq:class1}
\pi(\gamma_{t^5,t^4u})&=&\ell\\
\label{eq:class2} \pi(\gamma_{t^4u,t^3u^2})&=&\ell+3\gamma\\ 
\pi(\sigma_{t,u})&=&\ell+4\gamma \label{eq:class3}
\end{eqnarray}
where $t$ and $u$ are distinct elements of $\{v,w,x,y,z\}$.

\bigskip
Our interest is in defining and computing Gopakumar-Vafa invariants of classes
in $\mathbf{H}$ in order to compare to the calculations of \cite{HJKOV}.

\begin{lem}\label{lem:finite}
 Fix $\delta\in\mathbb{Z}^2$.  Then the set of curve classes $\beta\in M$ with $\pi(\beta)=\delta$ is finite.
\end{lem}

\bigskip\noindent
\emph{Proof.\/} We write $\delta=m\ell+n\gamma$ and may assume $m\ge0$, $n\ge0$, and $(m,n)\ne(0,0)$, otherwise the statement of the lemma
is trivial.  We can write 
\begin{equation}
\beta=\sum m_{\mathbf{t}}\ell^{\mathbf{t}}+\sum n_{i,j}\gamma_{\mathbf{m_1},\mathbf{m_2}}+\sum p_{\mathbf{s},\mathbf{t}}\sigma_{\mathbf{s},\mathbf{t}}
\label{eq:beta}
\end{equation} 
with $m_{\mathbf{t}},n_{i,j},p_{\mathbf{s},\mathbf{t}}>0$, suppressing from the notation the implicit set of 3 variables needed to define $\gamma_{\mathbf{m_1},\mathbf{m_2}}$.  In writing (\ref{eq:beta}), we have adopted the convention that we always use $\ell$ in place of a class of the form $\gamma_{t^5,t^4u}$.
Then $\pi(\beta)=\delta$ implies $\sum m_{\mathbf{t}}+\sum n_{i,j}+\sum  p_{\mathbf{s},\mathbf{t}}=m$ and $3\sum n_{i,j}+4\sum  p_{\mathbf{s},\mathbf{t}}=n$, which has only finitely many solutions for $m_{\mathbf{t}}, n_{i,j}$ and
$p_{\mathbf{s},\mathbf{t}}$.

Lemma \ref{lem:finite} ensures that the two parameter projection of \cite{HJKOV} makes mathematical sense: their invariants are simply a sum of finitely many Gopakumar-Vafa invariants.

\bigskip
For $\delta\in \mathbf{H}$ we define 
\begin{equation}\label{eq:projgw}
n^g_{m,n}=\sum_{\pi(\beta)=m\ell+n\gamma}n^g_\beta,
\end{equation}
which is a finite sum by Lemma~\ref{lem:finite}.  In
(\ref{eq:projgw}), $n^g_\beta$ is the genus $g$ GV invariant associated with 
the curve class $\beta\in H_2(X_\psi,\mathbb{Z})$. The GV invariants will be discussed in more detail in Section~\ref{sec:enum}.

\section{Enumerative Geometry}\label{sec:enum}
In this section, we calculate some genus~0 Gromov-Witten invariants.  In principle we can algorithmically compute the genus~0 Gromov-Witten invariants for the full 101 parameter model using the toric mirror theorem~\cite{givental}, but we leave that for future work.  It is conceivable that there may be a way to modify the toric mirror theorem to apply directly to our two-parameter
family.  That would certainly be worth doing.

\subsection{Generalities}
In lieu of calculating genus~0 Gromov-Witten invariants, we calculate
the genus~0 Gopakumar-Vafa invariants \cite{GV2}, which is easier in
examples.  Let $\beta\in H^2(X_\psi,\mathbb{Z})$.  Mathematically, we follow 
\cite{KatzGV0} by defining the
genus~0 Gopakumar-Vafa invariant $n^0_\beta$ as the Donaldson-Thomas invariant
of the moduli space of stable sheaves $F$ with $[F]=\beta$ (which is known to 
be independent of the choice of polarization used to define stability).  This
definition is consistent with the mathematical definition of higher genus 
Gopakumar-Vafa invariants $n^g_\beta$ given later in \cite{MT}.

The Gromov-Witten invariants $N^0_\beta$ are related to the $n^0_\beta$  by the
Aspinwall-Morrison formula \cite{tftrc,MR1421397}\footnote{A proof of the multiple cover formula has only been written down for isolated curves.  However, in our cases,  it can be shown by symplectic techniques that the $\bar\partial$ operator can be deformed so that there are only finitely many embedded pseudoholomorphic curves, all isolated, and the number is given by the associated GV invariant, with signed counts corresponding to orientation choices.  This is easy for primitive curve classes since the moduli spaces are smooth and projective, but with more care the case of $2\gamma$ below can also be handled.}
\begin{equation}\label{eq:multiplecover}
  N^0_\beta=\sum_{k|\beta}\frac{n^0_{\beta/k}}{k^3}
\end{equation}
so we content ourselves with the calculation of the $n^0_\beta$.   

Furthermore, in all of our cases, the stable sheaves are of the 
form $\mathcal{O}_C$ for rational curves $C$ and so the moduli spaces of sheaves is just the moduli space of curves. Even better,
these moduli spaces $M_\beta$ are smooth, and the Donaldson-Thomas invariant is simply
\begin{equation}\label{eq:euler}
n^0_\beta=\left(-1\right)^{\dim\left(M_\beta\right)}e\left(M_\beta\right),
\end{equation}
where $e\left(M_\beta\right)$  is the topological euler characteristic.

\subsection{Degrees 0 and 1}
We now record our results for Gopakumar-Vafa invariants in the classes $m\ell+n\gamma$, with $m=0,1$.  These are the classes which project to points or lines in the singular model $Y_\psi$.

\begin{eqnarray}
n^0_{\gamma}&=&300\\
n^0_{2\gamma}&=&-440\\
%n^0_{3\gamma}&=&780\\
n^0_\ell&=&15\\
n^0_{\ell+\gamma}&=&-60\\
n^0_{\ell+2\gamma}&=&155
\end{eqnarray}

We begin with the classes of the form $n\gamma$.  By Lemma~\ref{lem:deg0} and its proof, these curves can be constructed from the $\gamma$ curves and the fibers of the divisors of the $A_4$ resolution.

\bigskip\noindent
{\bf $\gamma$.} The curves of class $\gamma$ contained in each rational ruled surface are just the components of the reducible fibers.  The curves of class $\gamma$ contained in each $dP_3$ are the 6 curves depicted in Figure~\ref{fig:one-divisor}.  So we simply count all of these curves.  

Each 2-simplex depicted in Figure~\ref{fig:triangulation} corresponds to one of the 10 singular points. We count 30 curves of class $\gamma$ corresponding to each 2-simplex.  Thus $n^0_{\gamma}=10\cdot 30 =300$.

\bigskip\noindent
{\bf $2\gamma$.} Again inspecting Figure~\ref{fig:one-divisor}, we see that any intersecting conifiguration of two curves of type $\gamma$ must lie in exactly one of the $dP_3$s or the ruled surfaces lying over a singular curve.  Since each curve $\gamma$ has self-intersection $-1$ in the surface just described, the union $\gamma_1+\gamma_2$ is a rational curve of self-intersection $0$ in a rational surface, hence it moves in a pencil.  The corresponding contribution to the GV invariant is $(-1)^1e(\mathbb{P}^1)=-2$
by (\ref{eq:euler}).

There are 3 such pairs of intersecting $-1$ curves for each of the 6 $dP_3$'s over each of the 10 singular points.

There are 4 ruled surfaces over each to the 10 singular curves.

Thus, $n^0_{2\gamma}=(-2)(10\cdot3\cdot6+ 10\cdot 4)=-440$.

\bigskip

We now turn to the curves of the form $\ell+n\gamma$. Since all of these project to lines in $Y_\psi$, we start by identifiying the lines in $Y_\psi$.

\begin{lem}\label{lem:line}
Any line in $Y_\psi$ is contained in one of the five $\mathbb{P}^2$'s defined by the simultaneous vanishing of $u$
and one of the coordinates $\{v, w, x, y, z\}$, as well as
eq.~\eqref{eq:first}.
\end{lem}

\medskip\noindent
\emph{Proof.\/} The proof is inspired by an argument in \cite{AK}.  Let $\ell\subset Y_\psi$ be a line and suppose that $\ell$ is not contained in any of these $\mathbb{P}^2$'s.  If $\ell$ is not contained in the hyperplane $u=0$, then it intersects that hyperplane at  a point $p\in \ell$.  It follows from (\ref{eq:second}) that each of the five hypersurfaces $\mathbf{t}=0$ must also contain $p$.  This is a contradiction, since it would imply that 6 hyperplanes with empty intersection all contain $p$.

Since $u$ is identically 0 on $\ell$, it follows from (\ref{eq:second}) that at least one of the $\mathbf{t}$ must vanish identically on 
$\ell$, which proves the lemma.

\bigskip
Now let $C$ be any irreducible curve of degree 1, and we claim that its class is in the cone generated by the classes of the curves $\ell_{\mathbf{t}}$, $\sigma_{\mathbf{s},\mathbf{t}}$, and  $\gamma_{\mathbf{m_1},\mathbf{m_2}}$. In this case $\overline{C}=\rho(C)\subset Y_\psi$ is a line. It follows that $C$ must either be contained in one of the exceptional divisors forming an $A_4$ resolution or one of the divisors $D_{\mathbf{s}^5}$.  We have already shown the claim in the first case.  In the second case, $[C]=[\ell_{\mathbf{s}}]$.

\bigskip\noindent
{\bf $\ell$.}  The classes $\beta$ which contribute to $n^0_\ell$  project to a line in $Y_\psi$, which must lie in at least one of the five $\mathbb{P}^2$'s described in Lemma~\ref{lem:line}.  The divisor $D_{x^5}$ is the proper transform of one such $\mathbb{P}^2$, and the lines in such a $\mathbb{P}^2$ have class $\ell$.  The lines are parametrized by a dual $\mathbb{P}^2$ and so contribute $+3$.
There are five such $\mathbb{P}^2$, giving $5\cdot 3=15$.  

We claim that there are no other contributions.  The only remaining possibility is for the projection to lie in more than one $\mathbb{P}^2$ and thus be a singular line.  In this case one component of the curve must be a section of one of the 4 ruled surfaces over the singular line.  But from (\ref{eq:class1}) and (\ref{eq:class2}) we infer that the Mori cone of $D_{x^4y}$ is generated by $\ell$ and $\gamma$, while the only curve of class $\ell$ in $D_{x^4y}$ is $\gamma_{x^5,x^4y}$, which we have already considered as a line in
$D_{x^5}\simeq \mathbb{P}^2$.  Similarly, the Mori cone of $D_{x^3y^2}$ is generated by $\ell+3\gamma$ and $\gamma$ by
(\ref{eq:class2}) and (\ref{eq:class3}), so 
$D_{x^3y^2}$ contains no curves of class $\ell$.  Hence, the final answer is $n^0_\ell=15$.

\bigskip\noindent
{\bf $\ell+\gamma$.}  As above, we must have a component which is a line in one of the 5 $\mathbb{P}^2$'s, and we look for all ways  to attach a curve $\gamma$.    A typically way is to take the line in $D_{x^5}$ passing through the point $D_{x^5}\cap D_{x^4y}\cap D_{x^4z}$ (there is a pencil of such lines) and attach $\gamma_{x^4y,x^4z}$.  All other cases are obtained by permutations of the underlying data.  Such curves are associated with a partial flag $p\in
\mathbb{P}^2$, where $p$ is one of the 10 singular points and $\mathbb{P}^2$ is one of the 3 $\mathbb{P}^2$'s containing a given singular point.   This
gives $n^0_{\ell+\gamma}=10\cdot3\cdot(-2)=-60$.
\

\bigskip\noindent
{\bf $\ell+2\gamma$.} Again, we must have a component which is a line in one of the 5 $\mathbb{P}^2$'s.  We look for all ways to attach two curves of class $\gamma$ or one curve of class 2$\gamma$.   For the first case, we choose one of the five $\mathbb{P}^2$'s which we identify with some $D_{\mathbf{t}^5}$ and a pair of the 6 singular points contained in that $\mathbb{P}^2$.  This specifies a line in $D_{\mathbf{t}^5}$ containing the two singular points.  Then we glue the curves of type $\gamma$ which meet each of the singular points.  This contributes $5\cdot15=75$.

We have already seen how the fibers of the ruled surfaces have class $2\gamma$.  A typical example is to glue a line $\ell\subset D_{x^5}$ containing a point $p\in \gamma_{x^5,x^4y}$ to a fiber $\varphi_{x^4y}$.  All other cases are given by permutation.  The data
$p\in\ell$ is parametrized by the blowup of the dual $\mathbb{P}^2$ at the point $\gamma_{x^5,x^4y}$ of this dual $\mathbb{P}^2$.  This blowup has euler characteristic $4$ and contributes $+4$ to the GV invariant.  
We have such a pencil for each choice of one of the 5 $\mathbb{P}^2$'s and one of the 4 singular lines contained in each $\mathbb{P}^2$.  This contributes $5\cdot 4\cdot 4=80$.

Combining these cases, we get $n^0_{\ell+2\gamma}=75+80=155$.

\subsection{Degree 2}

Here we only have partial results, but we note that they already agree with parallel calculations in \cite{HJKOV}.  We compute the invariants of those curves of class $2\ell+n\gamma$ which lie inside the union of the toric divisors $D_{\mathbf{m}}$, and denote the contribution of these curves by $n^{0,{\rm toric}}_{2\ell+n\gamma}$.  However, since we do not have an analogue of Lemma~\ref{lem:line} at our disposal, it is possible that $n^{0,{\rm toric}}_{2\ell+n\gamma}\ne n^{0}_{2\ell+n\gamma}$.

We record our results.
\begin{eqnarray}
n^{0,{\rm toric}}_{2\ell}&=&-30\\
n^{0,{\rm toric}}_{2\ell+\gamma}&=&150\\
n^{0,{\rm toric}}_{2\ell+2\gamma}&=&-500
\end{eqnarray}

\bigskip\noindent
{\bf $2\ell$.}  In this class, we clearly have the conics in any of the five toric divisors $D_{\mathrm{s}^5}$, each of which is a $\mathbb{P}^2$.  Since conics are parametrized by $\mathbb{P}^5$, the contribution is $5(-6)=-30$.

The argument that we used in the case $\beta=\ell$ shows that there are no other curves in $\cup D_{\mathbf{m}}$ representing $\beta=2\ell$.

\bigskip\noindent
{\bf $2\ell+\gamma$.}  Analogous to the case of $\ell+\gamma$, our curve must be a connected union of one of the conics of class $2\ell$ and a curve $\gamma$, to which is associated to one of 30 possible partial flags $p\in \mathbb{P}^2$.  The linear system of conics in $\mathbb{P}^2$ containing a point $p$ is a  $\mathbb{P}^4$.  So the invariant is $30\cdot 5=150$.

\bigskip\noindent
{\bf $2\ell+2\gamma$.}  As in the situation of $\ell+2\gamma$, we must have a conic $C\subset D_{\mathbf{s}^5}$ as a component.  There are then two possibilities.  Either $C$ is glued to two distinct $\gamma$ curves, each meeting $D_{\mathbf{s}^5}$ in one of the 6 singular points that it contains, or $C$ is glued to a fiber of a ruled surface $D_{\mathbf{s^4t}}$.

In the first case, we choose one of the five $\mathbb{P}^2$ and one of the 15 pairs of the 6 singular points in that $\mathbb{P}^2$ at which we glue a $\gamma$ curve.  The linear system of plane conics through a pair of points is a $\mathbb{P}^3$.  So the net contribution is $5\cdot 15\cdot (-4)=-300$.

In the second case, for each of the 5 divisors $D_{\mathbf{s}^5}$, we pick one of the 4 singular lines $\gamma_{\mathbf{s}^5,\mathbf{s}^4\mathbf{t}}$ contained in it. We get a curve of class $2\ell+2\gamma$ by choosing a point $p\in \gamma_{\mathbf{s}^5,\mathbf{s}^4\mathbf{t}}$, taking a conic $C\subset D_{\mathbf{s}^5}$ containing $p$, and gluing it to the fiber $\phi_{\mathbf{s}^4\mathbf{t}}$ containing $p$.  This moduli space is a $\mathbb{P}^4$ bundle over $\mathbb{P}^1$, contributing $-10$.  So the total contribution is $5\cdot 4\cdot(-10)=-200$.

Combining these two cases, we get $n^{0,{\rm toric}}_{2\ell+2\gamma}=-300-200=-500$.

\section{Five Parameter Projection}

In this section, we prove Conjecture~\ref{conj:main} in a five-parameter projection.  Let us first explain what this means.  

It follows from Batyrev's calculation \cite{Batyrev1993} (which used
toric constructions to represent curve classes, as we have done)
 that the curves  $\ell_{\mathbf{t}}$, $\sigma_{\mathbf{s},\mathbf{t}}$, and  $\gamma_{\mathbf{m_1},\mathbf{m_2}}$ identified in Conjecture~\ref{conj:main} span $H_2(X_\psi,\mathbb{Q})$.  Furthermore, the curves  $\sigma_{\mathbf{s},\mathbf{t}}$ are not needed by (\ref{eq:sigma}).

We let $V$ be the $\mathbb{Q}$-vector space with generators $\ell_t, \sigma_{s,t}$, and $\gamma_{m_1,m_2}$.  Then there is a natural surjective map $r:V \to H_2(X_\psi,\mathbb{Q})$. The kernel of $r$ is spanned by the relations (\ref{eq:relgamma})--(\ref{eq:sigma}) together with the analogous relations obtained by substitutions of the variables, described in Section~\ref{sec:curvesdivisors}.

The group $S_5$ acts on $Y_\psi$ by permuting the coordinates, and this
action lifts naturally to an action on $X_\psi$.  
Let $W$ be the $\mathbb{Q}$-vector space generated by the set of $S_5$-orbits of the $\ell_t, \sigma_{s,t}$, and $\gamma_{m_1,m_2}$.   There is a natural map $n:V \to W$.  Since the set of relations spanning $\ker(r)$ just described is invariant under $S_5$, we deduce a natural map
\begin{equation}
\rho: H_2(X_\psi,Q) \to \mathbf{N} := W/\left(n\left(\ker r\right)\right).
\end{equation}

%We let
%$\mathbf{N}=H_2(X_\psi,\mathbb{Z})^S_5$ be the $S_5$ invariants, and we
%let $\rho:H_2(X_\psi,\mathbb{Z})\to \mathbf{N}$ be the natural projection.  
Let $\mathbf{M}\subset \mathbf{N}\otimes\mathbf{R}$ be the cone spanned by the classes $\rho(\ell_{\mathbf{t}})$, $\rho(\sigma_{\mathbf{s},\mathbf{t}})$, and $\rho(\gamma_{\mathbf{m_1},\mathbf{m_2}})$.  Recall the Mori cone $M\subset H_2(X_\psi,\mathbb{R})$.

\begin{lem}
$\mathrm{dim}(\mathbf{N})=5$.
\end{lem}

\medskip\noindent
\emph{Proof.\/} As noted above, we can ignore the $\rho(\sigma_{\mathbf{s},\mathbf{t}})$.  Clearly, all of the curves $\ell_{\mathbf{s}}$ are identified by the $S_5$ action.  There are six $S_5$ orbits of curves $\gamma_{\mathbf{m_1},\mathbf{m_2}}$:

\begin{equation}
\gamma_{\mathbf{s^4t},\mathbf{s^4u}},
\gamma_{\mathbf{s^4t},\mathbf{s^3tu}},
\gamma_{\mathbf{s^3t^2},\mathbf{s^3tu}},
\gamma_{\mathbf{s^3t^2},\mathbf{s^2t^2u}},
\gamma_{\mathbf{s^3tu},\mathbf{s^2t^2u}},
\gamma_{\mathbf{s^2tu^2},\mathbf{s^2t^2u}}.
\end{equation}

However, the third relation in (\ref{eq:relgamma}) gives $\rho(\gamma_{\mathbf{s^4t},\mathbf{s^3tu}})=\rho(\gamma_{\mathbf{s^3tu},\mathbf{s^2t^2u}})$ since the first and last terms cancel in the projection, and similarly the fifth relation in (\ref{eq:relgamma}) gives $\gamma_{\mathbf{s^3t^2},\mathbf{s^2t^2u}}=\gamma_{\mathbf{s^2tu^2},\mathbf{s^2t^2u}}$, since the second and third terms cancel in the projection.  Thus the $\rho(\gamma_{\mathbf{m_1},\mathbf{m_2}})$ span a 4-dimensional space.  Including the $\rho(\ell_{\mathbf{s}})$, we get $\mathrm{dim}(\mathbf{N})=5$.

\begin{prop}
$\rho(M)=\mathbf{M}$.
\end{prop}

\noindent
\emph{Proof.}  Let $\tau\subset \mathbf{M}$ be the simplicial cone spanned by the $\rho(\gamma_{\mathbf{m_1},\mathbf{m_2}})$ and $\rho(\ell_{\mathbf{s}})$.  We compute the dual cone
$\tau^\vee$, expressing the edges in terms of $S_5$-invariant divisors on $X_\psi$.

\medskip\noindent
{\bf Claim.} Each of these five $S_5$-invariant divisors can be chosen to be effective.

\smallskip
We defer the calculation and complete the proof.

 Let $C$ be an irreducible curve in the mirror quintic.  If $C$ is contained in a toric divisor, we already know that $[C]$ is contained in the semigroup spanned by $\ell_{\mathbf{t}}$, $\sigma_{\mathbf{s},\mathbf{t}}$, and  $\gamma_{\mathbf{m_1},\mathbf{m_2}}$  by our previous calculation of the Mori cone of these toric surfaces.  So $\rho([C])\in\mathbf{M}$.

If $C$ is not contained in any toric divisor, then its intersection with each toric divisor is nonnegative, and in particular its intersection with each of the five $S_5$-invariant divisors is nonnegative.  But these intersection numbers are simply the coefficients of $\rho(C)$ as it is expressed as a linear combination of the edges of $\tau$.   Therefore $\rho(C)\in\tau\subset \mathbf{M}$.

It remains to compute the edges of $\tau^\vee$.    We can express the five generators of $\tau$ as 
\begin{equation}
\rho(\ell_{\mathbf{s}}),\rho( \gamma_{\mathbf{s^4t},\mathbf{s^4u}}),\rho(\gamma_{\mathbf{s^4t},\mathbf{s^3tu}}),\rho(\gamma_{\mathbf{s^3t^2},\mathbf{s^3tu}}),\rho(\gamma_{\mathbf{s^3t^2},\mathbf{s^2t^2u}})
\end{equation}

We have already constructed an effective toric divisor $D^x$ with $D^x\cdots\ell_{\mathbf{s}}=1$ and having intersection number zero with each $\gamma$ curve.  Therefore, symmetrizing $D^x$ gives an $S_5$-invariant divisor which spans the edge of $\tau^\vee$ dual to the face of $\tau$ which does not contain $
\rho(\ell_{\mathbf{s}})$.

\section{Conclusions}

We have studied the curves on the quintic mirror threefold, using toric
constructions to obtain curves, and using toric surfaces to see the
relations among them.  Those same toric surfaces allowed us to isolate
some of these curves as potentially extremal:  they are extremal on the
surface that contains them, but not (yet) known to be extremal on the
threefold itself.  We hope that this approach to studying curves will 
be useful in other contexts as well.

\section*{Acknowledgements}
S.K. thanks the Mathematical Sciences Research Institute
and D.R.M. thanks the Aspen Center for Physics for hospitality during
portions of this research.
We thank Kantaro Ohmori and especially Cumrun Vafa for useful discussions.
The work of S.K.\ was partially supported by National Science Foundation
Grants DMS-1502170 ad DMS-1802242 as well as DMS-1440140 at MSRI.  
The work of D.R.M.\ was partially supported by Simons Foundation Award 
\#488629, as part of the Simons Collaboration on Special Holonomy in
Geometry, Analysis, and Physics, as well as National Science Foundation
Grant PHY-1607611 at the Aspen Center for Physics.

%%  The bibliography

%%  If your bibliography is in BibTeX format, use the following setup:
%%  Style BST file for numbered citation:
%\bibliographystyle{plain}
%%  Bibliography file (usually `*.bib')
%\bibliography{pamq-bibliography}
%%
%%  or include bibliography directly:
%\begin{thebibliography}{9}
%%  Use \bibitem{r1} or \bibitem[Surname(2010)]{r1} (for authoryear case)
%%  Put author names in \textsc{} command in order to use small caps font
%
%\bibitem{}
%\textsc{}
%
%\end{thebibliography}

\bibliographystyle{amsplain-enspra}
\bibliography{egmq}

\end{document}